\documentclass[11pt,makeidx]{article}
\input amssym.def
\input amssym.tex
\usepackage{epsfig,psfig}
\def\biblitem#1{\bibitem{#1}}
\def\binom#1#2{{#1}\choose{#2}}

\def\slfrac#1#2{\hbox{\kern.1em %
 \raise.5ex\hbox{\the\scriptfont0 #1}\kern-.11em %
 /\kern-.15em\lower.25ex\hbox{\the\scriptfont0 #2}}}

\newcommand{\eqn}[1]{(\ref{#1})}
\newcommand{\hsp}{\hspace*{\parindent}}
\newcommand{\vsp}{\vspace{.1in}}

\newcommand{\eeq}{\end{equation}}
\newcommand{\beql}[1]{\begin{equation}\label{#1}}

\newcommand{\La}{\Lambda}
\newcommand{\Th}{\Theta}

\newcommand{\dT}{{\rm det}\,}

\newcommand{\sC}{{\cal C}}

\makeatletter
\def\@sect#1#2#3#4#5#6[#7]#8{\ifnum #2>\c@secnumdepth
     \def\@svsec{}\else
     \refstepcounter{#1}\edef\@svsec{\csname the#1\endcsname.\hskip .75em }\fi
     \@tempskipa #5\relax
      \ifdim \@tempskipa>\z@
        \begingroup #6\relax
          \@hangfrom{\hskip #3\relax\@svsec}{\interlinepenalty \@M #8\par}%
        \endgroup
       \csname #1mark\endcsname{#7}\addcontentsline
         {toc}{#1}{\ifnum #2>\c@secnumdepth \else
                      \protect\numberline{\csname the#1\endcsname}\fi
                    #7}\else
        \def\@svsechd{#6\hskip #3\@svsec #8\csname #1mark\endcsname
                      {#7}\addcontentsline
                           {toc}{#1}{\ifnum #2>\c@secnumdepth \else
                             \protect\numberline{\csname the#1\endcsname}\fi
                       #7}}\fi
     \@xsect{#5}}
\def\@begintheorem#1#2{\it \trivlist \item[\hskip \labelsep{\bf #1\ #2.}]}

\def\plain{plain}\ifx\fmtname\plain\csname fi\endcsname
     
     \input docstrip
     \preamble

     Do not distribute the stripped version of this file.
     The checksum in the header refers to the documented version.

     \endpreamble
     \generateFile{here.sty}{t}{\from{here.doc}{}}
     \endinput
\fi
\ifcat a\noexpand @\let\next\relax\else\def\next{%
    \documentstyle[here,doc]{article}\MakePercentIgnore}\fi\next
\ifx\@Hxfloat\@Hundef\else\expandafter\endinput\fi
\let\@Hxfloat\@xfloat
\def\@xfloat#1[{\@ifnextchar{H}{\@HHfloat{#1}[}{\@Hxfloat{#1}[}}
\def\@HHfloat#1[H]{%
\expandafter\let\csname end#1\endcsname\end@Hfloat
\vskip\intextsep\vbox\bgroup\def\@captype{#1}\parindent\z@
\ignorespaces}
\def\end@Hfloat{\egroup\vskip \intextsep}

\makeatother

\begin{document}
\begin{center}
{\Large {\bf A Note on Optimal Unimodular Lattices}} \\
\vspace{1.5\baselineskip}
{\em J. H. Conway} \\
\vspace*{.2\baselineskip}
Mathematics Department \\
Princeton University \\
Princeton, NJ 08544 USA \\
\vspace*{1\baselineskip}
{\em N. J. A. Sloane} \\
\vspace*{.2\baselineskip}
Information Sciences Research \\
 AT\&T Labs-Research \\
Florham Park, NJ 07932-0971, U.S.A. \\
\vspace{1.5\baselineskip}
Jan 27, 1998; last revised Aug 2, 1998 \\
\vspace{1.5\baselineskip}
{\bf ABSTRACT}
\vspace{.5\baselineskip}
\end{center}
\setlength{\baselineskip}{1.5\baselineskip}

The highest possible minimal norm of a unimodular lattice is determined
in dimensions $n \le 33$.
There are precisely five odd 32-dimensional lattices with the highest
possible minimal norm (compared with more than $8.10^{20}$ in dimension 33).
Unimodular lattices with no roots exist if and only if $n \ge 23$, $n \neq 25$.
\clearpage
\thispagestyle{empty}
\setcounter{page}{1}

\section{Introduction}
\hsp
The results stated in the abstract were announced in 1989
\cite{Me157}.
For their proof we define the {\em shadow} $S$
(cf. \cite{Me157}, \cite{Me158}) of an integral lattice
$\La$ as follows.
If $\La$ is odd, $S= ( \La_0 )^\ast \setminus \La^\ast$, where
the subscript ``0'' denotes ``even sublattice'' and ``$\ast$'' denotes
``dual'';
if $\La$ is even, $S= \La$.
It is immediate that the theta series $\Theta_S$ of the shadow of an $n$-dimensional odd lattice is related to theta series $\Th_\La$ of the lattice by
\beql{Eq1}
\Th_S (z) =
\sqrt{\dT \La} \left(
\frac{\eta}{\sqrt{z}} \right)^n \Th_\La \left( 1- \frac{1}{z} \right) ~,
\eeq
where $\eta = e^{\pi i /4}$.

If $\La$ is an odd unimodular lattice then we can write
\beql{Eq2}
\Th_\La (z) ~=~ \sum_{j=0}^{[n/8]}
a_j \Delta_8 (q)^j
\theta_3 (q)^{n-8j} ~,
\eeq
where $q= e^{\pi i z}$,
$$\Delta_8 (q) ~=~ q ~ \prod_{m=1}^\infty (1-q^{2m-1} )^8 (1-q^{4m})^8 ~,$$
and $\theta_2$, $\theta_3$, $\theta_4$ are the usual
Jacobi theta series
\cite[p. 102]{SPLAG}.
From \eqn{Eq1} and \eqn{Eq2} we have
\beql{Eq3}
\Th_S (z) = \sum_{j=0}^{[n/8]} \frac{(-1)^j}{16^j} a_j \theta_4 (q^2)^{8j}
\theta_2 (q)^{n-8j} =
\sum \beta_n q^n ~,~~~
{\rm (say)}.
\eeq

\section{Bounds}
\hsp
Suppose $\La$ is an odd unimodular lattice with the highest possible
minimal norm $\mu_n$.
Then $a_0 , \ldots, a_{\mu_n -1}$ are determined by \eqn{Eq2}.
Most of the bounds in Table 1 now follows from the following conditions:
$\Theta_\La$ and $\Theta_S$ must have nonnegative integer coefficients;
$\Theta_\La \equiv \Theta_S \equiv 1$ $(\bmod~2)$;
there is at most one nonzero $\beta_r$ for $r < (\mu_n +2)/2$;
$\beta_r =0$ for $r < \mu_n/4$;
and $\beta_r \le 2$ for $r < \mu_n /2$.
(The last three conditions follow from the fact if the four cosets
of $\La_0$ in $\La_0^\ast$ are $\La_0^{(0)}$,
$\La_0^{(1)}$,
$\La_0^{(2)}$,
$\La_0^{(3)}$, with $\La_0 = \La_0^{(0)}$,
$\La = \La_0^{(0)} \cup \La_0^{(2)}$, $S= \La_0^{(1)} \cup \La_0^{(3)}$,
then $u,v \in S$ implies $u \pm v \in \La$.)

We give two examples.
For dimension $n=9$, if minimal norm 2 were possible,
\eqn{Eq2} would imply $a_0 =1$, $a_1 = -18$, hence $\Theta_\La = 1+252q^2 + \ldots$, $\Theta_S = \frac{9}{4} q^{1/4} + \ldots$.
Since the coefficients of $\Th_S$ are not integers, we conclude that
$\mu_9 \le 1$.

For dimension $n=33$, if minimal norm 4 were possible,
\eqn{Eq2} would imply $a_0 = 1$, $a_1 = - 66$,
$a_2 = 660$, $a_3 = - 880$, hence $\Th_\La = 1+ (70290 +a_4)q^4 + \ldots$,
$\Th_S = (a_4/32768)q^{1/4} + (110 - 63a_4 / 32768)q^{9/4} + \ldots$.
From the initial term of $\Th_S$ we see that $a_4 =0$ or $2^{16}$,
but if $a_4 = 2^{16}$ the second term of $\Th_S$ is negative;
so $a_4 =0$ and $\Th_S = 110q^{9/4} + \ldots$.

The following additional argument is needed to eliminate this case
(and also for $n=13$).
At least one of $\La_0^{(1)}$ and $\La_0^{(3)}$ must contain 55 vectors of norm $9/4$.
Suppose $u,v \in \La_0^{(1)}$ with $u \cdot u = v \cdot v = 9/4$,
$u \neq \pm v$.
Then $v \cdot v \in \{ \pm 1/4, \pm 3/4 \}$.
But since $u \pm v \in \La$, $u-v \in \La_0$, the only possibility is
$u \cdot v = 1/4$.
The inner product matrix of these 55 vectors is therefore
$8I + J/4$, where $J$ is an all-1's
matrix, which has rank 55,
impossible in a 33-dimensional space.
Therefore $\mu_{33} \le 3$.

Table 1 gives $\mu_n$ for $n \le 40$ (for $34 \le n \le 39$ the value is either 3 or 4).
All the upper bounds for $8 \le n \le 40$, $n \neq 25$, follow from the above arguments.
In dimension 25 they give only $\mu_{25} \le 3$.
However, Borcherds has enumerated all 25-dimensional unimodular
lattices, showing that $\mu_{25} =2$
(\cite{Bor1}; the list is included in \cite{SN}).
\begin{table}[H]
\caption{Highest minimal norm $\mu_n$ of an $n$-dimensional unimodular lattice.}

$$
\begin{array}{ccccccccc}
n & 1 & 2 & 3 & 4 & 5 & 6 & 7 & 8 \\
\mu_n & 1 & 1 & 1 & 1 & 1 & 1 & 1 & 2 \\
~ \\
n & 9 & 10 & 11 & 12 & 13 & 14 & 15 & 16 \\
\mu_n & 1 & 1 & 1 & 2 & 1 & 2 & 2 & 2 \\
~ \\
n & 17 & 18 & 19 & 20 & 21 & 22 & 23 & 24 \\
\mu_n & 2 & 2 & 2 & 2 & 2 & 2 & 3 & 4 \\
~ \\
n & 25 & 26 & 27 & 28 & 29 & 30 & 31 & 32 \\
\mu_n & 2 & 3 & 3 & 3 & 3 & 3 & 3 & 4 \\
~ \\
n & 33 & 34 & 35 & 36 & 37 & 38 & 39 & 40 \\
\mu_n & 3 & \mbox{3-4} & \mbox{3-4} & \mbox{4} & \mbox{3-4} & \mbox{3-4} & \mbox{3-4} & 4 \\
\end{array}
$$
\end{table}

Lattices achieving the bounds in Table 1
are well-known for $n \le 24$, 32 and 40
\cite[Chap. 16]{SPLAG}.
Borcherds showed that there is a unique unimodular lattice in dimension 26 
with minimal norm 3,
and least one in dimension 27
(\cite{Bor1}; see also \cite[3rd ed.]{SPLAG}, \cite{SN}).
Bacher and Venkov \cite{BV27} showed that there are
exactly three such lattices in dimension 27 and exactly 38 in dimension 28.

In dimension 30 a unimodular lattice $L_{30}$ of minimal norm 3 may be obtained by
gluing together two copies of $\sqrt{2} L$, where $L$ is the unimodular
lattice $A_{15}^+$.
A 29-dimensional example can be obtained by taking an appropriate
vector $v_4 \in L_{30}$ with $v_4 \cdot v_4 = 4$, and
projecting $\{u \in L_{30} : u \cdot v_4 \equiv 0 ~ ( \bmod~2 ) \}$ onto
$v_4^\perp$.
Examples in dimensions 32 and 31 may be similarly obtained from
$\sqrt{2} D_{16}^+$.

The Conway-Thompson theorem \cite[p. 46]{MH} shows
that unimodular lattices with minimal norm $\ge 3$ exist in all dimensions
$n \ge 37$, and we have found explicit examples in dimensions 33-36 \cite{SN}.
A 36-dimensional unimodular lattice with minimal norm 4 was
recently constructed by G. Nebe (personal communication).
This expalins all the lower bounds in the table,
and establishes
the third result mentioned in the abstract.

The first gap in the table is at dimension 34,
where the theta series of a putative lattice with minimal norm 4 is
$$1+ 60180 q^4 + 2075904 q^5 + \ldots ~;
$$
its shadow would have theta series
$$204 q^{5/2} + 758200 q^{9/2} + 274625820 q^{13/2} + \ldots ~.$$
We do not know if such a lattice exists.

In \cite{Me157} it was also announced that \eqn{Eq2} and \eqn{Eq3} imply
that for all sufficiently large $n$, $\mu_n \le [(n+6)/10]$.
In fact a more delicate analysis yields $\mu_n \le 2[n/24] +2$, unless $n=23$
when $\mu_{23} \le 3$ --- see \cite{RS98}.

\section{Dimension 32}
\hsp
Suppose $\La$ is an odd unimodular lattice in dimension 32 having
minimal norm 4.
The argument of Section 2 shows that
\begin{eqnarray*}
\Th_\La (z) & = & 1 + 81344 q^4 + \ldots ~, \\
\Th_S (z) & = & 64q^2 + 144896 q^4 + \ldots ~.
\end{eqnarray*}
At least 32 vectors of $S$ must be in (say) the coset $\La_0^{(1)}$, but since their sum and difference is in $\La$,
such vectors must be orthogonal, and so both
$\La_0^{(1)}$ and $\La_0^{(3)}$ must contain exactly 32.
Suppose $u_1 = c (2,0, \ldots, 0), \ldots, u_{32} = c(0, \ldots, 0,2)$,
$c= 1/ \sqrt{8}$, are a set of 32 such vectors.
Then $\La$ contains the sublattice spanned by all vectors of the form
$c( \pm 4, \pm 4, 0^{30} )$.

We may now apply the argument of \cite[Chap. 12, p. 333]{SPLAG}, to construct a binary code $\sC$ of length 32, which must be doubly-even and self-dual.
There are precisely five such codes \cite{cp}, \cite{Me65}, so there are
therefore five possibilities for $\La$.
If $\sC$ is one of the five codes, the corresponding lattice $\La$
is constructed from it in the same way that the Leech lattice is constructed
from the Golay code:
$\La$ is spanned by the vectors
$c(-3, 1^{31} )$, $c(2u)$ for $u \in \sC$,
and $c( \pm 4^2, 0^{30} )$.
$\sC$ has 620 minimal vectors, of weight 8, so $\La$ contains
$2^7 .620 + 2^2 {\binom{32}{2}} = 81344$ minimal vectors, of norm 4.

Finally, we discuss odd unimodular 33-dimensional lattices.
The total mass $M= \sum | Aut \, \La |^{-1}$ of this genus is
$1.407 \ldots \times 10^{21}$ \cite[Chap. 16]{SPLAG}.
Let $M_{2k}$ be the mass of the lattices with exactly $2k$ vectors
of norm 1 or 2, $k = 0, 1,\ldots$~, so that $M= M_0 + M_2 + M_4 + M_6 +
\ldots$~.
Although for even $n$ the average theta series
\beql{Eq100}
\frac{1}{M} ~ \sum_{\La} \frac{\Th_\La (q)}{| Aut ~\La |}
\eeq
of the genus of odd unimodular lattices is given in
\cite[p. 386]{Cas}, \cite[Chap. 12]{Gro} and \cite[Section 7.4]{Ran},
no analogous formula seems to be known for odd $n$.
However, Eric Rains has pointed out that the results of \cite[p. 70]{Shim}
imply that for both even and odd $n ~>~ 4$ the average theta series is given by
\beql{Eq101}
\theta_3(z)^n \sum_{j = 0}^{[n/4]}~ c_j ~ (g_2(z)^j ~+~ h_2(z)^j)
\eeq
where
$$g_2(z) ~=~ 16 \, q \, \prod_{m = 1}^\infty ~
\left( \frac {1+q^{2m}}{1+q^{2m-1}} \right) ^8,
$$
$$h_2(z) ~=~ \prod_{m = 1}^\infty ~
\left( \frac {1-q^{2m-1}}{1+q^{2m-1}} \right) ^8,
$$
and the $c_j$ are chosen so that the coefficients $\alpha_i$
in the $q$-expansion of
$$
\theta_3(z)^n \sum_{j = 0}^{[n/4]}~ c_j ~ g_2(z)^j
$$
satisfy $\alpha_{4i} = 2^{n-2} \alpha_i$ for all $i$, and
the coefficient of $q^0$ in \eqn{Eq101} is $1$.

For $n = 33$ we find that \eqn{Eq100} is 
$$
1 ~+~ \frac{15535133760578}{505245773078238529}~q
~+~ \frac{719890853572979520}{505245773078238529}~q^2 ~+~ \ldots ~,
$$
$$ ~=~ 1 ~+~ 0.0000307\ldots ~q ~+~ 1.4248\ldots ~ q^2 ~+~ \ldots ~.
$$
\vsp
Therefore the total number of vectors of norms 1 and 2
in all lattices in the genus is $1.425 \ldots \times M =$
$2M_2 + 4M_4 + 6M_6 + \ldots~ \ge 2 (M-M_0)$.
Hence $M_0 \ge 0.404 \ldots \times 10^{21}$, and so the number of lattices
with minimal norm 3 (the highest possible value) is $\ge 8 \times 10^{20}$.

\paragraph{Acknowledgement.}
We thank Eric Rains for helpful comments.

\end{document}